\documentclass[11pt]{article}
\begin{document}
\def\doublespaced{\baselineskip=\normalbaselineskip\multiply\baselineskip
 by 150\divide\baselineskip by 100}
\doublespaced
\def\lsim{~{\rlap{\lower 3.5pt\hbox{$\mathchar\sim$}}\raise 1pt\hbox{$<$}}\,}
\def\gsim{~{\rlap{\lower 3.5pt\hbox{$\mathchar\sim$}}\raise 1pt\hbox{$>$}}\,}
\def\thisday{~September~10, 1998 ~and~ hep-th/yymmnnn~~}
\def\thisday{~\today ~and~ hep-pth/yymmnnn~~}


\vbox to 3.5cm {
\vfill
}
\begin{center}{\Large{Spatial Rotation of the Fractional Derivative in Two Dimensional Space}}\\
Ehab Malkawi  \\
{{Department of Physics,\\
United Arab Emirates University\\
 Al Ain, UAE}}\\
 emalkawi@uaeu.ac.ae
\end{center}
\vskip -10pt
\begin{center}
Abstract
\end{center}
The transformation of the partial fractional derivatives under spatial rotation in $R^2$ are derived for the Riemann-Liouville and Caputo
definitions. These transformation properties link the observation of physical quantities, expressed through fractional derivatives, with respect to different coordinate systems (observers). It is the hope that such understanding could shed light on the physical interpretation of fractional derivatives. Also it is necessary to able to construct interaction terms that are invariant with respect to equivalent observers.


\vskip 6pt



\vskip 3pt

\newpage

\section{Introduction}

Fractional calculus deals with differentiation and integration to arbitrary real or complex orders. The idea is as old as the integer-order calculus.  Extensive mathematical discussion of fractional calculus can be found in Refs.~[1]-[5] and references therein.
The techniques of fractional calculus have been applied to wide range of fields, such as physics, engineering, chemistry, biology, economics, control theory, signal image processing, groundwater problems, and many others.

Physical applications of fractional calculus span a wide range of topics and problems (for a review see Refs.~[6]-[9] and references therein). Generalizing fractional calculus to several variables, multidimensional space, and  generalization of  fractional vector calculus has been reported [10]-[13].  Also, progress has been reported on generalization of Lagrangian and Hamiltonian systems [14,15].

Despite the significant progress in applying fractional calculus to a wide range of physical problems, there is still a lack of satisfactory geometric and physical interpretation
of fractional calculus, in comparison with the
simple interpretations of their integer-order counterparts (see Ref.~[16] and references therein).
Effort has been devoted to relate fractional calculus and fractal geometry [17,18].
A different approach to geometric interpretation of fractional calculus is based on the idea of the contact of $\alpha-$th order [19]. However, a satisfactory interpretation is still missing [16].

A study of the symmetry of physical systems described by fractional calculus and investigating the transformation properties of fractional derivative operators
under specific groups can shed light on the nature of such operators. For example, studying the transformation of the fractional operators under point transformation (or change of variables) can be found in Ref.~[20, 21], even though the class of variables substitution preserving the form of the fractional derivative is found to be very narrow. Generalizing the definition of fractional derivative of a function with respect to another function can be used to study general transformation of point transformation [22].

The question that we are trying to address in this work is basically the following: Given a quantity that is described by a fractional derivative of some function in a specific coordinate system $(x,y)$, how can this quantity be expressed in a rotated coordinate system $(x^\prime,y^\prime)$? Thus, we try to relate the same fractional quantity in two equivalent coordinate systems. In this work, we only consider the cartesian space coordinates in two dimensions, $(x,\,y)$, and adhere to the  notion of time as invariant in all coordinate systems (nonrelativistic). In other words, we are looking into the transformation of fractional derivative operators and their effect on physical quantities under the rotation group $SO(2)$. We are not interested in constructing the irreducible representations of $SO(2)$, as done in the standard integer-order differentiation case [23], rather we focus on the orthonormal basis of the two dimesnional cartesian space of $R^2$. Intuitively, quantities expressed through fractional derivative operators are expected to behave differently from scalar, vector, and tensor quantities.

In this work  we investigate the transformation properties of the fractional derivative and its action on an invariant scalar field, under space rotation in two dimensions, and using both the Riemann-Liouville and Caputo definitions.
We compare between the Riemann-Liouville and Caputo definitions in relating the same physical quantity under rotation. We check if there is a major difference in the transformation properties between the two definitions.
In Ref.~[24], the transformation properties of the Riemann-Liouville fractional derivative of a scalar field under infinitesimal transformation of $SO(2)$ are derived. The current work is more general and can be applied to any function and is not bounded to infinitesimal transformation. Also, we include the case of the Caputo definition of fractional derivative. Their special result agrees with our general results.
In section 2, we give a brief introduction to fractional calculus to lay out our notation.
In section 3, we lay out the $SO(2)$ transformation properties of fields and their derivatives expressed in Cartesian coordinates. In section 4, we apply our method to the Riemann-Liouville and Caputo fractional derivatives and layout their corresponding transformations under $SO(2)$. Finally, in section 5 we give a brief discussion of our results.

\section{Fractional Calculus}

For mathematical properties of fractional derivatives and integrals one can consult Refs.[1]-[5] and the references therein.
In this section we
lay out the notation used in the rest of this work. We consider real analytic functions on $R^2$ and only discuss the Riemann-Liouville and Caputo definitions of the fractional derivative.
Let $f(x,y)$ to be a real analytic function in a specific domain in the Euclidian space $R^2$; $f: R^2 \rightarrow R$.
The $x$-partial fractional derivative of order $\alpha$ (keeping $y$ constant) and where $a$ is the
lower limit of $x$  is written as $_aD_x^\alpha f(x,y)$. Similarly for the $y$-partial fractional derivative of order $\alpha$ (keeping $x$ constat) and where $b$ is the lower limit of $y$ is written as $_bD_y^\alpha f(x,y)$. Since the $f(x,y)$ is analytic then the partial fractional derivatives are assumed to commute, i.e., $\left[ _aD_x^\alpha,\, _bD_y^\alpha\right]=0$.\\
\textbf{Definition 2.1.} \emph{The Cauchy's repeated integration
formula of the $n$th-order integration of the function $f(x,y)$ along $x$, keeping $y$ constant, can be written as
\begin{eqnarray}
_aD_x^{-n} f(x,y)&=& \int_a^x dx_{n-1} \int_a^{x_{n-1}} dx_{n-2}\dots \int_a^{x_{1}} f(x_0,y)\,dx_0
\\ \nonumber
&=&\frac{1}{(n-1)!}\int_a^x \frac{f(u,y)\,du}{{(x-u)}^{1-n}}\,\, .
\end{eqnarray}}
\emph{A similar relation for the\, $n$th-order integration of the function $f(x,y)$ along $y$, keeping $x$ constant,
\begin{eqnarray}
_bD_y^{-n} f(x,y)&=v& \int_b^y dx_{n-1} \int_b^{y_{n-1}} dy_{n-2}\dots \int_b^{y_{1}} f(x,y_0)\,dy_0
\\ \nonumber
&=&\frac{1}{(n-1)!}\int_b^y \frac{f(x,u)\,du}{{(y-u)}^{1-n}}\,\, .
\end{eqnarray}}
\textbf{Definition 2.2.} \emph{The Riemann-Liouville fractional integration of order $\alpha<0$ and along $x$, keeping $y$ constant,
is defined as
\begin{equation}
_aD_x^\alpha f(x,y)= \frac{1}{\Gamma(-\alpha)} \int_a^x \frac{f(u,y)\,du}{{(x-u)}^{1+\alpha}} \,\, ,
\label{eq0}
\end{equation}
and along $y$, keeping $x$ constant, is
\begin{equation}
_bD_y^\alpha f(x,y)= \frac{1}{\Gamma(-\alpha)} \int_b^y \frac{f(x,u)\,du }{{(y-u)}^{1+\alpha}}\, \, ,
\label{eq0}
\end{equation}
where $\Gamma(.)$ is the Gamma function.}

We consider $a=b=0$ and drop them from the notation henceforth.
\textbf{Definition 2.3.} \emph{The Riemann-Liouville partial fractional derivatives of the order $\alpha>0$, where $n-1<\alpha<n$ and $n\in N$, are defined as
\begin{eqnarray}
D_x^\alpha f(x,y) = D_x^n\, {D_x^{\alpha-n}} f(x,y) =
\frac{1}{\Gamma(n-\alpha)} \int_0^x \frac{f(u,y)\, du}{{(x-u)}^{\alpha-n+1}}\,\, , \\
D_y^\alpha f(x,y) = D_y^n\, {D_y^{\alpha-n}} f(x,y) =\frac{1}{\Gamma(n-\alpha)} \int_0^y \frac{f(x,u)\, du}{{(y-u)}^{\alpha-n+1}}\,\, ,
\end{eqnarray}
where to unify notation we used $D_x^n\equiv \frac{d^n}{dx^n}$ and $ D_y^n\equiv \frac{d^n}{dx^n}$.}

The leibniz rule applied to the Riemann-Liouville partial fractional derivatives can be written as [1]-[5],
\begin{eqnarray}
D_x^\alpha \left(f(x,y)\, g(x,y)\right) = \sum_{k=0}^{\infty} \frac{\Gamma(\alpha+1)}{\Gamma(k+1)\Gamma(\alpha-k+1)} D_x^{\alpha-k}f(x,y)
\, D_x^k g(x,y)\, ,\label{lib1} \\
D_y^\alpha \left(f(x,y)\, g(x,y)\right) = \sum_{k=0}^{\infty} \frac{\Gamma(\alpha+1)}{\Gamma(k+1)\Gamma(\alpha-k+1)} D_y^{\alpha-k}f(x,y)
\,D_y^k g(x,y)\, .
\label{lib2}
\end{eqnarray}
\textbf{Definition 2.4. }
\emph{The Caputo partial fractional derivatives of order $\alpha>0$, where $n-1<\alpha<n$ and $n\in N$, are defined as
\begin{eqnarray}
  ^{C}D_x^\alpha f(x,y)   \equiv  \frac{1}{\Gamma(n-\alpha)} \int_0^x   \frac{ D_u^n f(u,y)\, du}{{(x-u)}^{\alpha-n+1}} \,\, ,
  \label{cap1}\\
  ^{C}D_y^\alpha f(x,y)   \equiv    \frac{1}{\Gamma(n-\alpha)} \int_0^y \frac{ D_u^n f(x,u)\, du}{{(y-u)}^{\alpha-n+1}}\,\, .
\label{cap2}
\end{eqnarray}}
The Riemann-Liouville and Caputo definitions of the fractional derivative are related [1]-[5]
\begin{eqnarray}
  ^{C}D_x^\alpha f(x,y)   =   D_x^\alpha f(x,y) -\sum_{k=0}^{n-1}\frac{x^{k-\alpha}}{\Gamma(k+1-\alpha)} D_x^k f(x,y) {\Large|}_{x=0}  \,\, ,
  \label{recap1}\\
  ^{C}D_y^\alpha f(x,y)   =   D_y^\alpha f(x,y) -\sum_{k=0}^{n-1}\frac{y^{k-\alpha}}{\Gamma(k+1-\alpha)} D_y^k f(x,y) {\Large|}_{y=0} \,\, .
\label{recap2}
\end{eqnarray}

As an example, for the function $f(x,y)= x^\beta y^\lambda$, where $\beta, \lambda >-1$ we find
\begin{eqnarray}
D_x^\alpha  \left(x^\beta \, y^\lambda\right)&=& \frac{\Gamma(\beta+1)}{\Gamma(\beta-\alpha+1)}\left( x^{\beta-\alpha}\, y^{\lambda}
\right)\, , \\ \nonumber
 D_y^\alpha \left(x^\beta \, y^\lambda\right) &=& \frac{\Gamma(\lambda+1)}{\Gamma(\lambda-\alpha+1)}\left( x^{\beta}\, y^{\lambda-\alpha}
\right)\, , \\ \nonumber
D_x^\alpha D_y^\alpha \left(x^\beta \, y^\lambda\right) &=& \frac{\Gamma(\beta+1)\Gamma(\lambda+1)}{\Gamma(\beta-\alpha+1)\Gamma(\lambda-\alpha+1)}\left( x^{\beta-\alpha}\, y^{\lambda-\alpha}
\right)\, .
\end{eqnarray}
The Caputo partial fractional derivatives give the same result for $x^\beta \, y^\lambda$,
except for $\beta,\lambda\in N  <n$, where in this case the Caputo partial fractional derivatives vanish, as expected from Eqs.~(\ref{cap1},\,\ref{cap2}).

\section {Transformation of Functions and Integer-Order Derivatives Under Spatial Rotation}

We consider the Euclidian space $R^2$ and work with orthonormal basis of the Cartesian coordinate system.
The distinction between covariant and contravariant quantities disappears and we drop
this notation in our discussion.
As discussed before, we are not interested in discussing the irreducible representations of $SO(2)$ [23], instead we
keep working with the Cartesian coordinates $(x,y)$ since we want to express the partial fractional derivatives in the same way.
Consider the two Cartesian coordinate systems $xy$ and $x^\prime y^\prime$
sharing the same origin $(0,0)$.
The $(x^\prime\, , y^\prime)$ coordinates are related
to the $(x\,,y)$ coordinates though a counter-clockwise rotation by an angle $\phi$
about their identical axis of rotation $z^\prime=z$. For $\phi=0$ the two
systems of coordinates are overlapping. In this work we consider the case of passive
transformation of the coordinate system, thus we consider the observation of the same
quantity from two different coordinate frames, without actually rotating physical quantities.
The two coordinate systems are related by the transformation matrix
\begin{eqnarray}
R(\phi)= \left[ {%
\begin{array}{cc}
\cos\phi & \sin\phi \\
-\sin\phi & \cos\phi%
\end{array}
} \right] \, .
\end{eqnarray}
For the active transformation case the angle $\phi$ would be replaced by $(-\phi)$.
The coordinates $(x^{\prime }\, , y^{\prime })$ are related to the coordinates $(x,y)$
through
\begin{eqnarray}
\left[ {%
\begin{array}{c}
x^{\prime }  \\
y^{\prime }
\end{array}%
}\right] =\left[ {%
\begin{array}{cc}
\cos \phi  & \sin \phi  \\
-\sin \phi  & \cos \phi
\end{array}%
}\right] \left[ {%
\begin{array}{c}
x  \\
y
\end{array}%
}\right]\, .
\label{trans1}
\end{eqnarray}

We present some of the properties of the functional transformations of the Lie group as applied to function $f(x,y)$.
 The generator of the $SO(2)$ rotation (angular momentum), around the $z$ axis, $L$ is written as,
\begin{eqnarray}
L=y D^1_{x}-x D^1_{y} \, .
\end{eqnarray}
The group elements of the rotation group are written, in the exponential form, as $ e^{\phi L}$, where the exponential is defined as an infinite sum of differential operators, explicitly
\begin{eqnarray}
 e^{\phi L}\equiv \sum^\infty_{n=0}\frac{\phi ^{n}}{n!} \,L^{n}\, ,
\end{eqnarray}
where $\phi$ is the angle of rotation.
Since $L\, x= y$ and $ L\, y=-x$, it is easy to derive the action of the element $e^{\phi L}$ on the coordinates $x$ and $y$, explicitly
\[
e^{\phi L}\,x=\cos \phi \,x+\sin \phi \,y =x^{\prime }\, ,
\]%
\[
e^{\phi L}\,y= -\sin \phi \,x+\cos \phi \,y =y^{\prime }\,.
\]%
Thus we retrieve the standard result in Eq.~(\ref{trans1})\thinspace\ .\\
\textbf{Remark 3.1. }
In general, for $n,m \in Z^+$,
\begin{eqnarray}
{x^{\prime }}^{n}\,{y^{\prime }}^{m} = \left(e^{\phi L}\,x^n\right) \left(e^{\phi L}\,y^m\right) =e^{\phi L}\,\left(x^{n}\,y^{m}\right)\, .
\label{eq19}
\end{eqnarray}
 This result emerges from proving the following two Lemmas.\\
\textbf{Lemma I}: Given any two functions $A(x,y)$ and $B(x,y)$ that are analytic (of Class $C^\infty$) in their common domain, then
 \begin{eqnarray}
 L^n\left( A \,B\right) =\sum_{k=0}^{n}\frac{n!}{k!\, (n-k)!} \left(L^k\, A\right) \left(L^{n-k}\, B\right)\, ,
 \end{eqnarray}
 where $n\in Z^+$.
This is just the general Leibniz rule applied to the differential operator $L$ and the proof is similar using mathematical induction.
 The second important Lemma  is\\
\textbf{Lemma II}: Given any two functions $A(x,y)$ and $B(x,y)$ that are analytic (of Class $C^\infty$) in their common domain, then
\begin{eqnarray}
e^{\phi L} \left ( A(x,y) \, B(x,y)\right ) = \left(e^{\phi L} \, A(x,y) \right )\, \left(e^{\phi L} \,B(x,y)\right)\, .
\label{eq22}
\end{eqnarray}
\emph{Proof.}  We first write
\begin{eqnarray}
e^{\phi L} \left ( A(x,y) \, B(x,y)\right ) =\sum^{\infty}_{n=0} \frac{\phi^n}{n!}\, L^n\, \left ( A(x,y) \, B(x,y)\right )\\ \nonumber
 =\sum_{n=0} \frac{\phi^n}{n!} {\sum_{k=0}}^{n}\frac{n!}{k!\, (n-k)!} \left(L^k\, A\right) \left(L^{n-k}\, B\right)\, .
\end{eqnarray}
Noting that $1/(n-k)!$ vanishes for $k>n$ we rewrite the above results as
\begin{eqnarray}
e^{\phi L} \left ( A(x,y) \, B(x,y)\right )= \sum^{\infty}_{n=0} \sum^{\infty}_{k=0} \frac{\phi^n}{k!\, (n-k)!} \left(L^k\, A\right) \left(L^{n-k}\, B\right)\, .
\end{eqnarray}
Shifting the order of the sum and then shifting the index $n\rightarrow n+k$, we write
\begin{eqnarray}
e^{\phi L} \left ( A(x,y) \, B(x,y)\right )= \sum^{\infty}_{k=0} \sum^{\infty}_{n=-k} \frac{\phi^{n+k}}{k!\, n!} \left(L^k\, A\right) \left(L^{n}\, B\right)\, .
\end{eqnarray}
Noting that $1/n!$ vanishes for $n=-k,-k+1,\dots, -1$, we reach the final conclusion
\begin{eqnarray}
e^{\phi L} \left ( A(x,y) \, B(x,y)\right )= \sum^{\infty}_{k=0} \frac{\phi^{k}}{k!} \left(L^k\, A\right)
\sum^{\infty}_{n=0}  \frac{\phi^{n}}{n!} \left(L^{n}\, B\right)\\ \nonumber
=\left(e^{\phi L} \, A(x,y) \right )\, \left(e^{\phi L} \,B(x,y)\right)\, .
\end{eqnarray}
Thus, Eq.~(\ref{eq19}) follows.

Based on Eq.~(\ref{eq19}) we can discuss the transformation of functions under rotation. A function $\Psi (x,y)$ that is analytic at $(0,0)$,  can be always expanded as
\[
\Psi (x,y)=\sum_{n,m}\frac{a_{n\,m}}{n!\,m!}\, \,x^{n}\, y^{m} \, .
\]%
In the $x^{\prime }\,y^{\prime }$ coordinates frame, the transformed function $\Psi ^{\prime
}(x^{\prime },y^{\prime })$ can also be written as
\[
\Psi ^{\prime }(x^{\prime },y^{\prime })=\sum_{n,m}\frac{a_{nm}^\prime }{n!\,m!} \,
{x^{\prime}}^{n}\,{y^{\prime }}^{m}
=e^{\phi L}\, \sum_{n,m}\frac{a_{nm}^{\prime
}}{n!m!}\,\,\,{x}^{n}\,{y}^{m}\,.
\]%
If the constants $a_{nm}$ are independent of the coordinates, then we
develop the following important result.\\
\textbf{Remark 3.2. }
A function $\Psi (x,y)$ that is analytic in a specific domain, where its dependence on the rotation angle $\phi$ comes only through its
arguments $(x,y)$, is transformed as
\begin{eqnarray}
\Psi ^{\prime }(x^{\prime },y^{\prime })=e^{\phi L}\,\Psi (x,y)\, .
\end{eqnarray}

A scalar function $\Psi (x,y)$ that is invariant under spatial rotation, must transform as $\Psi^{\prime }(x^{\prime },y^{\prime })=\Psi
(x,y)$, thus
\begin{eqnarray}
L\, \Psi (x,y)=xD^1_{y}\Psi -yD^1_{x}\Psi =0\,,
\label{eq15}
\end{eqnarray}
assuming the dependence on the rotation angle $\phi$ occur only through the argument $(x,y)$.

Next we investigate the transformation of derivatives, one notes the following
\begin{eqnarray}
D^1_{x^\prime} \,{x^\prime}^n &=& n\, {x^\prime}^{n-1}= n\,e^{\phi\, L} x^{n-1} =e^{\phi L}\,
D^1_{x}\, x^n \\ \nonumber
&=& e^{\phi L} \,D^1_{x}\, e^{-\phi L}\,e^{\phi
L}\, x^n= e^{\phi L} \,D^1_{x}\, e^{-\phi L}\, {x^\prime}^n \, .
\end{eqnarray}
Thus we derive the important result
\begin{eqnarray}
D^1_{x^\prime} = e^{\phi L}\, D^1_x \, e^{-\phi L}\, ,
\end{eqnarray}
and similarly for
\begin{eqnarray}
D^1_{y^\prime} = e^{\phi L}\, D^1_y \, e^{-\phi L}\, .
\end{eqnarray}
The following commutation relations are easy to derive
\[
\left[ D^1_{x},L\right] =-D^1_{y}\,\, ,\,\, \,
\left[ D^1_{y},L\right] =D^1_{x}\, .
\]

Transformation of higher derivatives can be derived, for example
\begin{eqnarray}
D^2_{x^\prime} = e^{\phi L}\, D^1_{x}\, e^{-\phi L} \,
e^{\phi L}\, D^1_{x}\, e^{-\phi L} = e^{\phi L}\,
D^2_{x}\, e^{-\phi L}\, .
\end{eqnarray}
In general we conclude\\
\textbf{Remark 3.3. } The partial derivatives transform under spatial rotation as
\begin{eqnarray}
& &D^n_{x^\prime} = e^{\phi L} \, D^n_{x}\, e^{-\phi L}\,\, , \,\,\,\,
D^n_{y^\prime} = e^{\phi L}\, D^n_{y}\, e^{-\phi L}\, , \\ \nonumber
& & D^n_{x^\prime} D^m_{y^\prime} = e^{\phi L} \, D^n_{x}
\,D^m_{y}\, e^{-\phi L}\, ,
\end{eqnarray}
where $n,m \in N$.

As an illustration, one can show that
\begin{eqnarray}
\left [D^2_{x},L\right ]=-2 D^1_{x}{D^1_y} \, , \,\,\,\,
\left [D^2_{y},L\right ]=+2 D^1_{x} {D^1_y}\, .
\end{eqnarray}
Thus,
\begin{eqnarray}
\left [D^2_{x}+ D^2_{y},L\right ]=0\, ,
\end{eqnarray}
indicating that the Laplace operator $D^2_{x}+ D^2_{y}$
is a scalar operator, as expected.

An important remark is to check the transformation of the operator $L$ itself under spatial rotation.
\begin{eqnarray}
L^\prime = e^{\phi L}\, L\, e^{-\phi L} = L\, .
\end{eqnarray}
Thus the operator $L$ is a scalar itself, a result that can be verified easily. This is an important point as it allows for the definition of inverse transformation between the $xy$ and $x^\prime y^\prime$ coordinate frames.

\section{Transformation of the Fractional Derivative Under Spatial Rotation}

In this section we generalize the results of the previous section to fractional calculus, in particular to the Riemann-Liouville and Caputo definitions of the fractional derivatives. We consider only real functions and real values of the fractional order $\alpha>0$. This restriction limits the allowed domain space $(x,\, y)$. This is clear if we consider the transformation of $x^\beta$ and $y^\beta$, where $\beta\in R$.
\[
e^{\phi L}\,x^\beta={\left(\cos \phi \,x+\sin \phi \,y \right)}^\beta =\, {x^\prime}^\beta ,
\]%
\[
e^{\phi L}\,y^\beta = {\left(-\sin \phi \,x+\cos \phi \,y\right)}^\beta ={y^\prime}^\beta\,.
\]%
Therefore, we require $x,y>0$ and $x^\prime, y^\prime>0$, which is equivalent to only consider the allowed domain: $\left[ (x,y):  x,y>0 \cap y/x > \tan\phi \right]$. Despite the restriction on the domain, one can still study the transformation of the fractional derivatives and deduce their properties. An extension to the whole domain can be pursued in the complex plane once a clear understanding of such transformation is acquired in the restricted domain.

Based on the results of the previous section, we have\\
\textbf{Remark 4.1. }
In the domain $\left[x, y >0 \cap y/x > \tan\phi\right]$, the quantity $x^{\beta}\,y^{\lambda}$ transforms as
\begin{eqnarray}
{x^{\prime }}^{\beta}\,{y^{\prime }}^{\lambda} = \left(e^{\phi L}\,x^\beta\right) \left(e^{\phi L}\,y^\lambda\right) =e^{\phi L}\,\left(x^{\beta}\,y^{\lambda}\right)\, ,
\label{trans2}
\end{eqnarray}
where $\beta,\lambda\in R$.
Next we look at the transformation of the partial fractional derivatives. Using the Riemann-Liouville and Caputo definition of the fractional derivative we note the following
\begin{eqnarray}
{{D}^\alpha_{x^\prime}} \,{x^\prime}^\beta &=& \frac{\Gamma(\beta+1)}{\Gamma(\beta-\alpha+1)} \, {x^\prime}^{\beta-\alpha}
= \frac{\Gamma(\beta+1)}{\Gamma(\beta-\alpha+1)}\,e^{\phi\, L} x^{\beta-\alpha}
\\ \nonumber
 &=& e^{\phi L}\,
{{D}^\alpha_{x}} \, x^\beta = e^{\phi L} \,{{D}^\alpha_{x}} \, e^{-\phi L}\,e^{\phi
L}\, x^\beta= e^{\phi L} \,{{D}^\alpha_{x}} \, e^{-\phi L}\, {x^\prime}^\beta \, .
\end{eqnarray}
\textbf{Remark 4.2.}
The partial fractional derivatives according to the Riemann-Liouiville and Caputo definitions transform as
\begin{eqnarray}
{{D}^\alpha_{x^\prime}} &=& e^{\phi L}\, {{D}^\alpha_{x}} \, e^{-\phi L}\, , \\ \nonumber
{{D}^\alpha_{y^\prime}} &=& e^{\phi L}\, {{D}^\alpha_{y}} \, e^{-\phi L}\, .
\end{eqnarray}

Considering an infinitesimal angle transformation, $\phi << 1$, and keeping only linear terms of $\phi$, we write
\begin{eqnarray}
{{D}^\alpha_{x^\prime}} \approx {{D}^\alpha_{x}} + \phi \left[L, {{D}^\alpha_{x}}\right] \, , \\ \nonumber
{{D}^\alpha_{y^\prime}} \approx {{D}^\alpha_{y}} + \phi \left[L, {{D}^\alpha_{y}}\right] \, .
\end{eqnarray}
\textbf{Remark 4.3. }
Given an analytic function $f(x,y)$, its partial fractional derivatives transform as
\begin{eqnarray}
{{D}^\alpha_{x^\prime}} f^\prime(x^\prime,y^\prime)= e^{\phi L} {{D}^\alpha_{x}} f(x,y)\, ,\\ \nonumber
{{D}^\alpha_{y^\prime}} f^\prime(x^\prime,y^\prime)= e^{\phi L} {{D}^\alpha_{y}} f(x,y)\, .
\end{eqnarray}
Considering the infinitesimal angle approximation, $\phi<<1$, and keeping only linear terms of $\phi$, we write
\begin{eqnarray}
 {{D}^\alpha_{x^\prime}} f^\prime(x^\prime,y^\prime)\approx  D^\alpha_x f(x,y) + \phi L {{D}^\alpha_{x}} f(x,y)\, ,\label{trans4} \\
 {{D}^\alpha_{y^\prime}} f^\prime(x^\prime,y^\prime)\approx  D^\alpha_y f(x,y) + \phi L {{D}^\alpha_{x}} f(x,y)\, .
\end{eqnarray}
Thus we derived the transformation of the partial fractional derivatives according to the Riemann-Liouville and Caputo definitions.
Next we consider the transformation of the partial fractional derivatives as applied to a scalar field $\Psi(x,y)$ and using both Riemann-Liouville and Caputo definitions.

\subsection{ Transformation of the Riemann-Liouville Partial Fractional Derivatives of a Scalar Field }
 Consider a scalar function $\Psi(x,y)$ that is invariant under spatial rotation, i.e.,
$L\Psi(x,y)=0$, then according to Eq.~(\ref{trans4})
\begin{eqnarray}
 {{D}^\alpha_{x^\prime}} \Psi^\prime(x^\prime,y^\prime)\approx  D^\alpha_x \Psi(x,y) + \phi L {{D}^\alpha_{x}} \Psi(x,y)\, ,
\end{eqnarray}
where
\begin{eqnarray}
 L{{D}^\alpha_{x}} \Psi^\prime(x^\prime,y^\prime)= y {D^1_x} {{D}^\alpha_{x}} \Psi(x,y)  -x {D^1_y} {{D}^\alpha_{x}}  \Psi(x,y)\, .
\end{eqnarray}
Since $y$ and $D^1_y$ commute with $D^\alpha_x$ we can write the above equation as
\begin{eqnarray}
 L{{D}^\alpha_{x}} \Psi^\prime(x^\prime,y^\prime) = {D^1_x} {{D}^\alpha_{x}} \left(y\Psi\right)-x {{D}^\alpha_{x}} D^1_y \Psi\, .
\end{eqnarray}
Rewriting the above equation using commutation relations, we find
\begin{eqnarray}
 L{{D}^\alpha_{x}} \Psi^\prime(x^\prime,y^\prime) = \left[{D^1_x}, {{D}^\alpha_{x}} \right]\left(y\Psi\right)
 -\left[x , {{D}^\alpha_{x}}\right] \left(D^1_y \Psi\right) + {{D}^\alpha_{x}} \left(y D^1_x -x D^1_y\right)\Psi\, .
\end{eqnarray}
The last term vanishes since $L\Psi=0$. Thus we arrive at
\begin{eqnarray}
 L{{D}^\alpha_{x}} \Psi(x,y) = \left[{D^1_x}, {{D}^\alpha_{x}} \right]\left(y\Psi\right)
 -\left[x , {{D}^\alpha_{x}}\right] \left(D^1_y \Psi\right)\, .
 \end{eqnarray}
Note that from the properties of the fractional derivative, $D^m_x D^\alpha_x = D_x^{m+\alpha}$, while $\left[{D^1_x}, {{D}^\alpha_{x}} \right]\neq 0$ except for $\alpha\in Z^+$ . We use the known relation [1]-[5]
\begin{eqnarray}
D^\alpha_x D^1_x \Psi(x,y) = D_x^{\alpha+1} \Psi(x,y) - \frac{ x^{-\alpha-1}\left({ \Psi(0,y)}\right)}{\Gamma(-\alpha)}\, .
\end{eqnarray}
Thus, we conclude
 \begin{eqnarray}
 \left[{D^1_x}, {{D}^\alpha_{x}} \right] \left(y\,\Psi(x,y)\right) =  + \frac{y\,\Psi(0,y) x^{-\alpha-1}}{\Gamma(-\alpha)}\,\, .
 \label{comm1}
 \end{eqnarray}
 Using Leibniz rule in Eq.~(\ref{lib1}), substituting $f(x)= \Psi(x,y)$ and  $g(x)=x$, one can show that
 \begin{eqnarray}
 \left[x, {{D}^\alpha_{x}} \right] D^1_y\Psi(x,y) = -\alpha D_x^{\alpha-1}\,D^1_y\Psi(x,y)\, .
 \label{comm2}
 \end{eqnarray}
 Thus we conclude that
 \begin{eqnarray}
 {{D}^\alpha_{x^\prime}} \Psi^\prime(x^\prime,y^\prime)\approx  D^\alpha_x \Psi(x,y) + \phi \left( \alpha {{D}^{\alpha-1}_{x}} D^1_y\Psi +
 \frac{y\, \Psi(0,y)\, x^{-\alpha-1}}{\Gamma(-\alpha)}   \right)\, .
\label{trans5}
\end{eqnarray}
A similar relation holds for $y$
\begin{eqnarray}
 {{D}^\alpha_{y^\prime}} \Psi^\prime(x^\prime,y^\prime)\approx  D^\alpha_y \Psi(x,y) - \phi \left( \alpha {{D}^{\alpha-1}_{y}} D^1_x\Psi +
 \frac{x\,\Psi(x,0)\, y^{-\alpha-1}}{\Gamma(-\alpha)}   \right)\, .
\label{trans6}
\end{eqnarray}
Note that for $\alpha \in Z^+$, the last term in the above result vanishes. Also for $0<\alpha<1$, the term $D^{\alpha-1}_{y}$ represents
a fractional integral \\
\textbf{Remark 4.1.1}
For $\alpha=1$ we recover the standard transformation of the gradient operator,
\begin{eqnarray}
{{D}^1_{x^\prime}} \Psi^\prime(x^\prime,y^\prime) = D^1_x \Psi(x,y) + \phi \, D^1_y\Psi(x,y)\, ,\\ \nonumber
{{D}^1_{y^\prime}} \Psi^\prime(x^\prime,y^\prime) = D^1_y \Psi(x,y) - \phi \, D^1_x\Psi(x,y)\, . \nonumber
\end{eqnarray}
\textbf{Remark 4.1.2}
For $\alpha=2$ we recover the known scalar invariance of the Laplacian operator
\begin{eqnarray}
\left ( D^2_{x^\prime} +D^2_{y^\prime}\right) \Psi^\prime(x^\prime,y^\prime) = \left ( D^2_{x} +D^2_{y}\right) \Psi(x,y)\, .
\end{eqnarray}

\subsection{ Transformation of the Caputo Partial Fractional Derivatives of a Scalar Field}

To derive the transformation of the Caputo partial fractional derivatives we consider $0<\alpha<1$, and use the relation between the Riemann-Liouville and Caputo definitions, as expressed in Eqs.~(\ref{recap1},\,\ref{recap2}). Also we use the results of Eqs.~(\ref{comm1},\,\ref{comm2}). One can show that
\begin{eqnarray}
\left[x, {^C}D_x^\alpha\right]\left(D^1_y\Psi\right) = \left[x,D^\alpha_x\right]\left(D^1_y\Psi\right) - \frac{x^{1-\alpha}\,D^1_y\Psi(0,y)}{\Gamma(1-\alpha)}\, ,
\end{eqnarray}
\begin{eqnarray}
\left[D^1_x,{^C}D^\alpha_x\right]\left( y \Psi\right) = \left[D^1_x,D^\alpha_x\right]\left(y \Psi\right)  - \frac{x^{-\alpha-1} y\,\Psi(0,y)}{\Gamma(-\alpha)}\, .
\end{eqnarray}
Therefore, the Caputo partial fractional derivative transform, under infinitesimal angle $\phi<<1$, as
\begin{eqnarray}
 {^C{D}^\alpha_{x^\prime}} \Psi^\prime(x^\prime,y^\prime)\approx  {^C}D_x^{\alpha} \Psi(x,y) + \phi \left( \alpha {{D}^{\alpha-1}_{x}} D^1_y\Psi +
 \frac{ x^{1-\alpha}\, D^1_y\Psi(0,y)}{\Gamma(1-\alpha)}   \right)\, .
 \label{trans7}
\end{eqnarray}
\begin{eqnarray}
 {{D}^\alpha_{y^\prime}} \Psi^\prime(x^\prime,y^\prime)\approx  {^C}D^\alpha_y \Psi(x,y) - \phi \left( \alpha {{D}^{\alpha-1}_{y}} D^1_x\Psi +
 \frac{ y^{1-\alpha}\, D^1_x\Psi(x,0)}{\Gamma(1-\alpha)}   \right)\, .
\end{eqnarray}
Note that $\alpha-1<0$ and thus $D^{\alpha-1}_{x}$ and $D^{\alpha-1}_{y}$ represent fractional integrals.

\section{Discussion and Conclusions}
In section 4, we derived the transformation of the Riemann-Liouville and Caputo partial fractional derivatives under spatial rotation. 
We also
derived the transformation of the Riemann-Liouville and Caputo partial fractional derivatives of a scalar function.
We observe no fundamental difference between the transformations of the Riemann-Liouville and Caputo partial fractional derivatives. The fact that the Caputo fractional derivative of a constant vanishes does not lead to a significant difference in its transformation. The difference we observe is that in the Riemann-Liouville transformation, there is a term that depends on the value of the function at the boundary. While, for the Caputo transformation, a similar term depends on the normal derivative of the function at the boundary.
We observe that the Riemann-Liouville and Caputo partial fractional derivatives of a scalar function do not transform as the components of a gradient vector, as in the case of the integer-order, $\alpha=1$, partial derivatives. Thus they are classified to have a new form of transformation, different from the scalar, vector, and tensor quantities. They represent new entities that deserve to be studied thoroughly in future work.  

Next we provide few examples of scalar functions.
Consider the constant scalar function $f(x,y)=1$, According to Eqs.~(\ref{trans5},\,\ref{trans6}), the Riemann-Liouville transformations of the fractional derivatives are, for infinitesimal $\phi<<1$,
\begin{eqnarray}
 D^\alpha_{x^\prime} (1) \approx  D^\alpha_x (1)  + \phi \left(
 \frac{y\,\, x^{-\alpha-1}}{\Gamma(-\alpha)}   \right)\, ,\\ \nonumber
 D^\alpha_{y^\prime} (1) \approx  D^\alpha_y (1)  - \phi \left(
 \frac{x\,\, y^{-\alpha-1}}{\Gamma(-\alpha)}   \right)\, .\nonumber
\end{eqnarray}
One can easily check that the result is as expected, for example
\begin{eqnarray}
D^\alpha_{x^\prime} (1)= \frac{{x^\prime}^{-\alpha}}{\Gamma(-\alpha)}\approx \frac{{\left(x+\phi \,y\right)}^{-\alpha}}{\Gamma(-\alpha)} \approx D^\alpha_{x} (1)+
\phi\frac{ y\,x^{-\alpha-1}}{\Gamma(-\alpha)}\, ,
\end{eqnarray}
as required. It is easy to show that the following quantities are invariant,
\begin{eqnarray}
& & x^\alpha D^\alpha_{x} (1)\, , \,\,\,\,\,
 y^\alpha D^\alpha_{y} (1) \\ \nonumber
& & y^{-\alpha} D^\alpha_{x} (1) - x^{-\alpha} D^\alpha_{y} (1)  \, , \,\,\,\,
x^{2+\alpha} D^\alpha_{x} (1) + y^{2+\alpha} D^\alpha_{y} (1)\, .
\end{eqnarray}
One can use $(dx)^\alpha$ and $(dy)^\alpha$ instead of $x^\alpha$ and $y^\alpha$ in the above invariant forms. 
The Caputo partial fractional derivatives of a constant field vanish.

Consider the scalar function $\Psi(x,y) = x^2+y^2$. The Riemann-Liouville and Caputo $x$-partial fractional derivatives are, respectively
\begin{eqnarray}
D^\alpha_x \Psi(x,y)&=& \frac{2}{\Gamma(3-\alpha)} x^{2-\alpha} + y^2 \frac{x^{-\alpha}}{\Gamma(1-\alpha)}\, .\\ \nonumber
^C{D^\alpha_x} \Psi(x,y) &=& \frac{2}{\Gamma(3-\alpha)} x^{2-\alpha}\, . \nonumber
\end{eqnarray}
Again it is easy to check that the transformations of the partial fractional derivatives, for infinitesimal angle $\phi<<1$, agree with Eqs.~(\ref{trans5},\,\ref{trans7})
 for both the Riemann-Liouville and Caputo fractional derivatives.
A similar transformation holds for $D^\alpha_y$ and $^C{D^\alpha_y}$ with $x$ and $y$ interchanged and $\phi \rightarrow -\phi$.
Once can easily check that the quantities
\begin{eqnarray}
& & x^{\alpha} D^\alpha_{y} \Psi(x,y) + y^{\alpha} D^\alpha_{x} \Psi(x,y)\, , \,\,\,\,\, x^{\alpha}y^\alpha D^\alpha_y \,D^\alpha_{y} \Psi(x,y) \,
\label{last1} ,\\
& & x^{\alpha}\,\, {^C}{D_y^\alpha}\Psi(x,y) + y^{\alpha}\,\, {^C}{D_x^\alpha}\Psi(x,y)\, .
\label{last2}
\end{eqnarray}
are invariant under spatial rotation. Again, one can use $(dx)^\alpha$ and $(dy)^\alpha$ instead of $x^\alpha$ and $y^\alpha$ in the above invariant forms.

For the function $\Psi(x,y) = (x^2+y^2)^2$, one can check the transformations correct. However, the quantities
in Eq.~(\ref{last1}) are not invariant separately. Nevertheless, it is possible to find a linear combination of the two quantities to form an invariant quantity, namely,
\begin{eqnarray}
x^{\alpha} D^\alpha_{y} \Psi(x,y) + y^{\alpha} D^\alpha_{x} \Psi(x,y)+A\left(x^{\alpha}y^\alpha D^\alpha_y \,D^\alpha_{y} \Psi(x,y) \right)\, .
\end{eqnarray}
where $A$ is some constant to make the combination invariant. Again, one can use $(dx)^\alpha$ and $(dy)^\alpha$ instead of $x^\alpha$ and $y^\alpha$ in the above invariant forms.

In Ref.~[24], the authors derived the infinitesimal transformation of the partial fractional derivatives of a scalar field, by using the function's Taylor expansion and then applying the Riemann-Liouville fractional derivative to each term of the expansion. Their result agrees with our general result.

\section*{Acknowledments}
The author would like to thank Prof. Yuri Luchko and Prof. Dumitru Baleanu for useful discussion.

\section*{References}

\begin{itemize}

\item[[1]]
  Fractional Integrals and Derivatives: Theory and Applications, S.G. Samko, A.A. Kilbas, O.I. Marichev, Gordon and Breach, New York (1993).

\item[[2]]
An introduction to the fractional calculus and fractional differential equations, K.S. Miller, B. Ross, John Wiley \& Sons, New York (1993).

\item[[3]]
The Fractional Calculus; Theory and Applications of Differentiation and Integration to Arbitrary Order, K B. Oldham, J. Spanier, Academic Press,  New York (1974).

\item[[4]]
Fractional Calculus: Integrations and Differentiations of Arbitrary Order, K. Nishimoto, University of New Haven Press, New Haven (1989).

\item[[5]]
 Fractional Differential Equations, I. Podlubny, Academic Press, San Diego (1999).
\item[[6]]
\label{Herrmann}
Fractional Calculus. An Introduction for Physicists, Richard Herrmann, World Scientific, Singapore  (2011).

\item[[7]]
 Applications of Fractional Calculus in Physics, R. Hilfer, World Scientific, Singapore (2000);

\item[[8]]
 Fractals and Fractional Calculus in Continuum Mechanics, A. Carpinteri, F. Mainardi (Eds.), Springer, Wien (1997);
 Physics of Fractal Operators, B. West, M. Bologna, P. Grigolini, Springer, New York (2003).
 
 \item[[9]]
Review of Some Promising Fractional Physical Models, V. E. Tarasov, International Journal of Modern Physics B, Volume 27, Issue 09, 2013;
 Hamiltonian Chaos and Fractional Dynamics, G.M. Zaslavsky, Oxford University Press, Oxford (2005).

\item[[10]]
Multi–dimensional solutions of space-time-fractional diffusion equations, A. Hanyga, Proc. R. Soc. Lond. 2002 vol. 458, 429-450;
The differentiability in the fractional calculus, F. Ben Adda, Nonlinear Anal., 47 (2001).

\item[[11]]
Fractional curl operator in electromagnetics, N. Engheta, Microwave and Optical Technology Letters, Volume 17, Issue 2, 1998;
Complex and higher order fractional curl operator in electromagnetics, Q.A. Naqvi, M. Abbas, Optics Communications, Vol. 241, 349-355, 2004.

\item[[12]]
Fractional vector calculus for fractional advection-dispersion, M.M. Meerschaert, J. Mortensen, S.W. Wheatcraft, Physica A, 367, pp. 181–190, (2006);
Fractional differential forms, K. Cottrill-Shepherd, M. Naber, J. Math. Phys., 42, pp. 2203-2212, (2001).

\item[[13]]
Fractional Generalization of Gradient Systems, V.E. Tarasov, Letters in Mathematical Physics, Volume 73, Issue 1, pp 49-58 (2005);
Applications of fractional exterior differential in three-dimensional space, Yong Chen, Zhen-ya Yan, Hong-qing Zhang, Appl. Math. Mechanics, 24 (3) (2003), pp. 256-260.

\item[[14]]
Nonconservative Lagrangian and Hamiltonian mechanics, F. Riewe, Physical Review E 53 (1996) 1890-1899 (1996);
Hamilton-Jacobi and fractional like action with time scaling, M.A.E. Herzallah, I.M. Sami, D. Baleanu, and M.R. Eqab,
Nonlinear Dynamics 66 549-555 (2011);
On the Fractional Hamilton and Lagrange Mechanics, A.K. Golmankhaneh, M.Y. Ali, and D. Baleanu, International Journal of Theoretical Physics 51, (2012).

\item[[15]]
Hamiltonian formulation of systems with linear velocities within Riemann–Liouville fractional derivatives, Sami I. Muslih,
D. Baleanu, Journal of Mathematical Analysis and Applications, Volume 304, Issue 2, (2005);
Lagrangian and Hamiltonian fractional sequential mechanics, M. Klimek, Czech. J. Phys., 52 (2002);
Formulation of Euler–Lagrange equations for fractional variational problems, O.P. Agrawal, J. Math. Anal. Appl., 272, (2002).

\item[[16]]
{lubni1}
Geometric and Physical Interpretation of Fractional Integration and Fractional Differentiation, Igor Podlubny
Fractional Calculus and Applied Analysis, Volume 5, Number 4 (2002).

\item[[17]]
Geometry of fractional spaces, G. Calcagni, Adv. Theor. Math. Phys. 16 (2012) 549-644;
 Physical and geometrical interpretation of fractional operators, M. Monsrefi-Torbati and J.K. Hammond, Journal of the Franklin Institute
Volume 335, Issue 6, (1998), Pages 1077-1086;
A fractional integral and its physical interpretation, R. R. Nigmatullin, Theoret. and Math. Phys., vol. 90, no. 3, 1992, pp. 242-251.

\item[[18]]
On physical interpretations of fractional integration and differentiation, R. S. Rutman, Theoret. and Math. Phys., vol. 105, no. 3, 1995, pp. 1509-1519;
A probabilistic interpretation of the fractional order differentiation, J.A. Tenreiro Machado, Fract. Calc. Appl. Anal., 8 (2003), pp. 73-80.

\item[[19]]
Geometric interpretation of the fractional derivative, F. Ben Adda, Journal of Fractional Calculus, vol. 11, May 1997, pp. 21-52.


\item[[20]]
Invariance of a Partial Differential Equation of Fractional Order under the Lie Group of Scaling Transformations, E. Buckwar,
Y. Luchko, Journal of Mathematical Analysis and Applications, Volume 227, Issue 1, 1 November 1998, Pages 81-97;
Invariance of a partial differential equation of fractional order under the
Lie group of scaling transformations, E. Buckwar, Yu. Luchko , J. Math. Anal. Appl., 1998, V. 227, P. 81-97.

\item[[21]]
Symmetry properties of fractional diffusion equations, R.K. Gazizov, A.A. Kasatkin, S.Yu. Lukashchuk, Physica Scripta. 2009. T136, 014016;
Group-Invariant Solutions of Fractional Differential Equations, R.K. Gazizov, A.A. Kasatkin, S.Yu. Lukashchuk, Nonlinear Science and Complexity,
(2011), pp. 51-59.

\item[[22]]
Fractional Differential Equations: Change of Variables and Nonlocal Symmetries, R.K. Gazizov, A.A. Kasatkin, S.YU. Lukashchuk,
Ufa  Mathematical Journal. Vol. 4. No 4 (2012). pp 54-67.

\item[[23]]
Cartesian Tensors, Jeffreys H. Cambridge University Press, (1931);
Mathematical Methods for Physicists, George B. Arfken, Hans J. Weber, Frank E. Harris,
Elsevier Academic Press, (2012);
Matrix Groups: An Introduction to Lie Group Theory, Baker Andrew, Springer (2003) ;
Differential Geometry, Lie Groups, and Symmetric Spaces (Graduate Studies in Mathematics), Sigurdur Helgason,
American Mathematical Society, (2001).


\item[[24]]
Transformation of Fractional Derivatives Under  Space Rotation. By E. Malkawi and A.A. Rousan. Published in International Journal of Applied Mathematics, Volume 16 No. 2, 175-185, 2004.

\end{itemize}

\end{document}